\documentclass[a4paper]{amsart}
\usepackage{graphicx}
\usepackage{amssymb}
\usepackage{amsmath}
\usepackage{amsthm}
\usepackage{amscd}
\usepackage[all,2cell]{xy}

\UseAllTwocells \SilentMatrices
\newtheorem{thm}{Theorem}[section]

\newtheorem{cor}[thm]{Corollary}
\newtheorem{lem}[thm]{Lemma}

\newtheorem{prop}[thm]{Proposition}
\theoremstyle{definition}

\theoremstyle{remark}
\newtheorem{rem}[thm]{\bf Remark}
\numberwithin{equation}{section}

\begin{document}
\title[Deriving Milnor's theorem on pullback rings]{Deriving Milnor's theorem on pullback rings}
\author[Xiao-Wu Chen, Jue Le] {Xiao-Wu Chen, Jue Le$^*$}

\date{\today}

\thanks{$^*$ The corresponding author}

\thanks{E-mail: xwchen$\symbol{64}$mail.ustc.edu.cn, juele$\symbol{64}$ustc.edu.cn}

\maketitle

\dedicatory{}%
\commby{}

\begin{abstract}
 The classical theorem of Milnor on pullback rings states that the category of projective modules over a pullback ring is equivalent to a certain category of gluing triples consisting of projective modules.  We prove an analogous result on the level of derived categories, where the equivalence has to be replaced by an epivalence.
\end{abstract}

\section{Introduction}
 Let $R_i$ and $R'$ be three rings with unit for $i=1,2$.  Let $\pi_i\colon R_i \rightarrow R'$ be two ring homomorphisms. The associated \emph{pullback ring}  is defined to be $R=\{(a, b)\in R_1\times R_2\; |\; \pi_1(a)=\pi_2(b)\}$, which is a subring of the product $R_1\times R_2$.

 Let us first recall the classical theorem in \cite{Mil} on  projective modules over the pullback ring $R$. For a ring $\Lambda$, we denote by $\Lambda\mbox{-Mod}$ the category of left $\Lambda$-modules, and by $\Lambda\mbox{-Proj}$ the full subcategory formed by projective modules. We denote by $\mathbf{D}(\Lambda\mbox{-Mod})$ the unbounded derived category of $\Lambda\mbox{-Mod}$.

  A \emph{gluing triple} $(X_1, X_2; c)$ consists of an $R_1$-module $X_1$, an $R_2$-modules $X_2$ and an $R'$-isomorphism  $c\colon R'\otimes_{R_1}X_1\rightarrow R'\otimes_{R_2} X_2$. This defines the category ${\rm Tr}_{\rm gl}(\pi_1, \pi_2)$ of gluing triples. For an $R$-module $M$, we have a natural gluing triple ${\rm Ind}(M)=(R_1\otimes_R M, R_2\otimes_R M; {\rm can}_M)$, where ${\rm can}_M$ is the canonical isomorphism between the induced $R'$-modules. This gives rise to the \emph{induction} functor
 $${\rm Ind}\colon R\mbox{-Mod}\longrightarrow {\rm Tr}_{\rm gl}(\pi_1, \pi_2).$$

 The following characterization theorem of projective modules over a pullback ring is fundamental; see \cite[Section~2]{Mil}.

 \vskip 5pt

 \noindent {\bf Theorem.} (Milnor) \emph{Assume that $\pi_1\colon R_1\rightarrow R'$ is surjective. Then the induction functor ${\rm Ind}$ restricts to an equivalence
 $$R\mbox{-{\rm Proj}} \stackrel{\sim}\longrightarrow \{(P_1, P_2; c)\in {\rm Tr}_{\rm gl}(\pi_1, \pi_2)\;| \; P_i\in R_i\mbox{-{\rm Proj}}\}.$$}

 The module categories of certain pullback rings are studied in \cite{NR, Le, AL}. In the representation theory of artin algebras, two important classes of algebras, namely string algebras and skew-gentle algebras, are realized as certain pullback algebras \cite{EHIS, BD}. Here, we notice that a radical embedding between algebras naturally gives rise to a pullback algebra. Milnor's theorem is useful in constructing new derived equivalences between pullback rings \cite{HX}.

Recall that a functor $F\colon \mathcal{C}\rightarrow \mathcal{D}$ is an \emph{epivalence},  provided that it is full and dense, and detects isomorphisms between objects. Consequently, it induces a bijection between the  sets of isomorphism classes of objects in $\mathcal{C}$ and $\mathcal{D}$.

For a radical embedding, it is shown that there are epivalences from the derived category of modules over the pullback ring to a certain category formed by gluing derived-triples of complexes; see \cite[Theorem~2.4]{BD04} and \cite[Theorem~4.2]{BD}. These epivalences are very useful to describe indecomposable objects in the derived category of the pullback ring.

The above mentioned results in \cite{BD04, BD} are analogous to Milnor's theorem. It is natural to expect that such an analogue might hold in a more general setting. We will confirm this expectation to some extent.

We define the category $\mathbf{D}{\rm Tr}_{\rm gl}(\pi_1, \pi_2)$ of \emph{gluing derived-triples} by replacing the module categories in defining ${\rm Tr}_{\rm gl}(\pi_1,\pi_2)$ with the corresponding derived categories. Using the derived tensor functors, we have the \emph{derived induction} functor
$${\rm Ind}^\mathbb{L}\colon \mathbf{D}(R\mbox{-Mod})\longrightarrow \mathbf{D}{\rm Tr}_{\rm gl}(\pi_1, \pi_2).$$

The following main result states an epivalence  in the derived setting; see Theorem~\ref{thm:main}. It is a derived analogue to Milnor's theorem.

\vskip 5pt

\noindent{\bf Theorem.}\quad \emph{Assume that $\pi_1\colon R_1\rightarrow R'$ is surjective and that $R'$ is finitely generated projective as the induced right $R_2$-module by $\pi_2$. Then the derived induction functor ${\rm Ind}^\mathbb{L}$ is full, whose kernel ideal is square zero. Consequently, there is an epivalence
$$\mathbf{D}(R\mbox{-{\rm Mod}})\longrightarrow {\rm Im}({\rm Ind}^\mathbb{L}),$$
where ${\rm Im}({\rm Ind}^\mathbb{L})$ denotes the essential image of ${\rm Ind}^\mathbb{L}$.}

\vskip 5pt

In general, it is not easy to characterize the essential image of ${\rm Ind}^\mathbb{L}$. When one considers right-bounded complexes, the corresponding essential image is more accessible. Then we  extend \cite[Theorem~2.4]{BD04} and a part of \cite[Theorem~4.2]{BD}. Here, we emphasize that the kernel ideal is square zero.

For the proof of the main result, we study the category of possibly non-gluing triples, and observe that it is equivalent to the module category over an upper triangular matrix ring $\Gamma$. Then the result follows from an explicit $1$-tilting $\Gamma$-module \cite{HR} arising from the pullback ring, and an epivalence  arising in the canonical recollement \cite{BBD} associated to $\Gamma$. We mention that the later epivalence is a special case of \cite[Theorem~A]{CL}; compare \cite[Appendix]{Gei}.

In comparison with the equivalence in Milnor's theorem, the epivalence in our main result is weaker. As mentioned above, the epivalence really appears in the canonical recollement associated to $\Gamma$. One reason for such an epivalence is that derived categories, or triangulated categories in general, are not rigid enough, and thus the middle category in a recollement is only determined by the outer ones up to an epivalence, not up to an equivalence; see \cite{CL}.

In conclusion, working with derived categories, an epivalence in our main result is reasonable. As is well known, to resolve the non-rigidity of a triangulated category, one might enhance them to stable $\infty$-categories. We mention that \cite[Theorem~28]{Tam} is quite analogous to our main result, which states an equivalence on the level of stable $\infty$-categories for a certain pullback of ring spectra.

The paper ia structured as follows. In Section~2, we recall Milnor's theorem. Under a certain additional condition, we extend it to an equivalence between the category of separated modules \cite{Le} over a pullback ring and the category of gluing triples; see Proposition~\ref{prop:ext-Mil}. In Section~3, we relate the category of triples to the category of modules over an upper triangular matrix ring $\Gamma$. Moreover, we obtain an explicit $1$-tilting $\Gamma$-module, which arises from the pullback diagram of rings; see Proposition~\ref{prop:tilt}.  In Section~4, we recall an epivalence arising from the canonical recollement associated to $\Gamma$,  and prove Theorem~\ref{thm:main}. The essential image of the derived induction functor is discussed in Proposition~\ref{prop:rb} and Remark~\ref{rem:image}.

\section{Modules over pullback rings}
In this section, we recall the induction and pullback functors. Under certain conditions, we extend
Milnor's theorem: these functors induce mutually inverse equivalences between the category of separated modules over a pullback ring and the category of gluing triples; see Proposition~\ref{prop:ext-Mil}.

Throughout this paper, we fix a pullback diagram of rings.
\begin{align}\label{equ:sq}
\xymatrix{
R\ar[d]_-{i_2} \ar[r]^-{i_1} & R_1\ar[d]^-{\pi_1}\\
R_2\ar[r]^-{\pi_2} & R'}
\end{align}
Therefore, given any $a_i\in R_i$ satisfying $\pi_1(a_1)=\pi_2(a_2)$, there exists a unique $a\in R$ satisfying $i_1(a)=a_1$ and $i_2(a)=a_2$.

For a ring $\Lambda$, we denote by $\Lambda\mbox{-Mod}$ the category of left $\Lambda$-modules, and by $\Lambda\mbox{-Proj}$ the full subcategory of projective modules. We denote by $\Lambda\mbox{-proj}$ the full subcategory formed by finitely generated projective modules.

By default, modules mean left modules. Therefore, right $\Lambda$-modules will be viewed as left modules over the opposite ring $\Lambda^{\rm op}$.

\subsection{The induction and pullback functors}

Following \cite{Mil}, we use triples to study $R$-modules. The category ${\rm Tr}(\pi_1, \pi_2)$ of triples is defined as follows: a \emph{triple} $(X_1, X_2; c)$ consists of left $R_i$-modules $X_i$ and an $R'$-morphism $c\colon R'\otimes_{R_1} X_1\rightarrow R'\otimes_{R_2} X_2$; a morphism $(f_1, f_2)\colon (X_1, X_2;c)\rightarrow (Y_1, Y_2; c')$ between triples is given by $R_i$-morphims $f_i\colon X_i\rightarrow Y_i$ satisfying
$$c'\circ (R'\otimes_{R_1} f_1)=(R'\otimes_{R_2} f_2)\circ c.$$
A triple $(X_1, X_2;c)$ is said to be \emph{gluing}, provided that the $R'$-morphism $c$ is an isomorphism. Denote by ${\rm Tr}_{\rm gl}(\pi_1, \pi_2)$ the full subcategory formed by gluing triples.

The following  \emph{induction} functor
$${\rm Ind}\colon R\mbox{-Mod}\longrightarrow {\rm Tr}(\pi_1, \pi_2)$$
sends an $R$-module $M$ to the canonical gluing triple $(R_1\otimes_{R} M, R_2\otimes_R M; {\rm can}_M)$, where
$${\rm can}_M\colon R'\otimes_{R_1}(R_1\otimes_{R} M)\longrightarrow R'\otimes_{R_2}(R_2\otimes_{R} M)$$
is the canonical $R'$-isomorphism. On the other hand, we have the \emph{pullback} functor
$${\rm Pb}\colon {\rm Tr}(\pi_1, \pi_2) \longrightarrow R\mbox{-Mod}$$
sending a triple $(X_1, X_2;c)$ to
$${\rm Pb}(X_1, X_2; c)=\{(x_1, x_2)\in X_1\times X_2\; |\; c(1\otimes x_1)=1\otimes x_2\}.$$
The $R$-module structure on ${\rm Pb}(X_1, X_2; c)$ is given by $a(x_1, x_2)=(i_1(a)x_1, i_2(a)x_2)$ for $a\in R$.

The following fact is well known.

\begin{lem}
The pair $({\rm Ind}, {\rm Pb})$ is adjoint.
\end{lem}

\begin{proof}
The natural isomorphism
$${\rm Hom}_R(M, {\rm Pb}(X_1, X_2;c)) \longrightarrow {\rm Hom}_{{\rm Tr}(\pi_1, \pi_2)}({\rm Ind}(M), (X_1, X_2;c))$$
sends $f$ to $(f_1, f_2)$. Here, the $R_i$-morphism $f_i\colon R_i\otimes_R M\rightarrow X_i$ is determined by $f_i(1\otimes m)={\rm pr}_i\circ f(m)$, where ${\rm pr}_i\colon {\rm Pb}(X_1, X_2;c)\rightarrow X_i$ denotes the projection on the corresponding entry.
\end{proof}

\begin{rem}
We describe the unit and counit of the above adjoint pair. For an $R$-module $M$, the unit
$$\eta_M\colon M\longrightarrow {\rm Pb}\circ {\rm Ind}(M)$$
sends $m$ to $(1\otimes m, 1\otimes m)$. Here, we view ${\rm Pb}\circ {\rm Ind}(M)$ as an $R$-submodule of $(R_1\otimes_R M)\oplus (R_2\otimes_R M)$. For a triple $(X_1, X_2; c)$, the counit
$$\epsilon_{(X_1, X_2;c)} \colon {\rm Ind}\circ {\rm Pb}(X_1, X_2; c)\longrightarrow (X_1, X_2;c)$$
is given by $R_i\otimes_R {\rm Pb}(X_1, X_2;c)\rightarrow X_i $, which sends $1\otimes (x_1, x_2)$ to $x_i$.
\end{rem}

The following fundamental result is contained in \cite[Section 2]{Mil}.

\begin{thm}\label{thm:Mil}
Assume that $\pi_1\colon R_1\rightarrow R'$ is surjective. Then the adjoint pair $({\rm Ind}, {\rm Pb})$ restricts to mutually inverse equivalences
$$R\mbox{-{\rm Proj}} \stackrel{\sim}\longrightarrow \{(P_1, P_2; c)\in {\rm Tr}_{\rm gl}(\pi_1, \pi_2)\;| \; P_i\in R_i\mbox{-{\rm Proj}}\}$$
and
$$R\mbox{-{\rm proj}} \stackrel{\sim}\longrightarrow \{(P_1, P_2; c)\in {\rm Tr}_{\rm gl}(\pi_1, \pi_2)\;| \; P_i\in R_i\mbox{-{\rm proj}}\}.$$
\end{thm}

\begin{rem}\label{rem:Mil}
(1) We observe that ${\rm Ind}(R)=(R_1, R_2; {\rm can}_R)$. By the adjoint pair $({\rm Ind}, {\rm Pb})$, we infer that the induction functor induces an isomorphism
$${\rm End}_R(R)\simeq {\rm End}_{{\rm Tr}(\pi_1, \pi_2)}({\rm  Ind}(R)),$$
both of which are isomorphic to $R^{\rm op}$. Moreover, the functor ${\rm Ind}$ commutes with arbitrary coproducts. It follows that ${\rm Ind}$ restricts a fully faithful functor
$${\rm Ind}\colon R\mbox{-{\rm Proj}} \longrightarrow \{(P_1, P_2; c)\in {\rm Tr}_{\rm gl}(\pi_1, \pi_2)\;| \; P_i\in R_i\mbox{-{\rm Proj}}\}.$$
Therefore, the essential part of Theorem~\ref{thm:Mil} claims that the surjectivity of $\pi_1$ implies the density of the above functors.

(2) We mention that the surjectivity condition of $\pi_1$ is necessary for Theorem~\ref{thm:Mil}. For an example, we let $k$ be a field, and assume that $k=R=R_1=R_2$ and $R'=k[t]/(t^2)$. The natural injective homomorphisms between them form a pullback diagram of rings. Consider the gluing triple $(k, k; c)$ satisfying $c(1\otimes 1)=(1+t)\otimes 1$. We observe  ${\rm Pb}(k, k;c )=0$. It follows that the triple $(k, k; c)$ is not induced from any $R$-module.
\end{rem}

\subsection{An extension of Theorem~\ref{thm:Mil}}

We will extend these equivalences under an additional condition on $\pi_1\colon R_1\rightarrow R'$. Recall from \cite{Le} that an $R$-module $M$ is \emph{separated}, provided that there exists an $R$-monomorphism  $M\rightarrow X_1\oplus X_2$ for some $R_i$-modules $X_i$.

The following result is implicitly contained in \cite[Corollary~3.3]{Le}; compare \cite[Proposition~1.1]{Wis}.

\begin{lem}\label{lem:Au} Let $M$ be an $R$-module. Then the following statements hold.
\begin{enumerate}
\item The unit $\eta_M$ is a monomorphism if and only if $M$ is separated.
\item Assume that $\pi_1$ is surjective. Then $\eta_M$ is an epimorphism.
\end{enumerate}
\end{lem}

\begin{proof}
We recall ${\rm Pb}\circ {\rm Ind}(M)\subseteq (R_1\otimes_R M)\oplus (R_2\otimes_R M)$ is an $R$-submodule.  We observe that any $R$-morphism $M\rightarrow X_1\oplus X_2$ factors through $\eta_M$. Then (1) follows immediately.

For (2), the surjectivity of $\pi_1$ gives rise to an exact sequence of $R$-bimodules
$$0\longrightarrow R \stackrel{\binom{i_1}{i_2}}\longrightarrow R_1\oplus R_2 \stackrel{(\pi_1, -\pi_2)}\longrightarrow R'\longrightarrow 0.$$
Applying $-\otimes_R M$ to it and identifying $R\otimes_R M$ with $M$, we obtain an exact sequence of $R$-modules
\begin{align}\label{equ:M}
M\longrightarrow (R_1\otimes_R M)\oplus (R_2\otimes_R M)\longrightarrow R'\otimes_R M\longrightarrow 0.
\end{align}
Since ${\rm Pb}\circ {\rm Ind}(M)$ coincides with the kernel of $(\pi_1, -\pi_2)\otimes_R M$, we infer that $\eta_M$ is always an epimorphism.
\end{proof}

We say that a two-sided ideal $I$ of a ring $\Lambda$ is \emph{universally superfluous}, provided that for any $\Lambda$-module $Z$, the $\Lambda$-submodule $IZ\subseteq Z$ is superfluous. For example, any nilpotent ideal is universal superfluous. If $\Lambda$ is left perfect, then any  ideal $I$  contained in ${\rm rad}(\Lambda)$, the Jacobson radical of $\Lambda$,  is universally superfluous; see \cite[Proposition~24.4(2)]{Lam}.

\begin{lem}\label{lem:Ba}
Let $(X_1, X_2; c)$ be a triple.
\begin{enumerate}
\item If $\epsilon_{(X_1, X_2;c)}$ is an isomorphism, then $(X_1, X_2;c)$ is gluing.
\item Assume that $\pi_1$ is surjective and that ${\rm Ker}\pi_1$ is a universally superfluous ideal of $R_1$. Then  $\epsilon_{(X_1, X_2;c)}$ is an isomorphism for any gluing triple $(X_1, X_2; c)$.
\end{enumerate}
\end{lem}

\begin{proof}
The statement in (1) is trivial. For (2), we set $I= {\rm Ker}i_2$ and $I_1={\rm Ker}\pi_1$.  We observe that $i_2\colon R\rightarrow R_2$ is also surjective and that $i_1$ restricts to an isomorphism $I\simeq I_1$. The $R$-module ${\rm Pb}(X_1, X_2;c)=M$ is given by the following pullback diagram
\[\xymatrix{
M\ar[d]_-{p_2} \ar[r]^-{p_1} & X_1\ar[d]^-{\xi}\\
X_2 \ar[r]^-{\kappa} & R'\otimes_{R_2} X_2,
}\]
where $p_i$ denotes the projection, $\kappa(x_2)=1\otimes x_2$ and $\xi(x_1)=c(1\otimes x_1)$. Therefore, $p_2\colon M\rightarrow X_2$ is surjective. It follows that $\kappa \circ p_2(M)=\xi\circ p_1(M)$ generates $R'\otimes_{R_2} X_2$ as an $R'$-module. Using the isomorphism $c$, we identify $R'\otimes_{R_2} X_2$ with $$R'\otimes_RX_1\simeq R_1/{I_1}\otimes_{R_1}X_1\simeq X_1/{I_1X_1}.$$
 Then we infer that $R_1p_1(M)+I_1X_1=X_1$. Since $I_1$ is universally superfluous, we have $R_1p_1(M)=X_1$. Consequently, we have
\begin{align}\label{equ:IM}
p_1(IM)=I_1p_1(M)=I_1R_1p_1(M)=I_1X_1.
\end{align}

 Since $I$ acts trivially on $X_2$, we have $IM\subseteq {\rm Ker}p_2$. On the other hand, $p_1$ restricts to an isomorphism ${\rm Ker}p_2\simeq {\rm Ker}\xi=I_1X_1$. In view of (\ref{equ:IM}), we infer that $IM={\rm Ker}p_2$. Consequently, the canonical map
 $$R_2\otimes_R M\longrightarrow  X_2,\quad 1\otimes(x_1, x_2)\mapsto x_2$$
  is an isomorphism, since both modules are isomorphic to $M/IM$. Therefore, the following $R'$-morphism
 $$\delta\colon R'\otimes_R M\longrightarrow R'\otimes_{R_2} X_2, \quad 1\otimes (x_1, x_2)\mapsto c(1\otimes x_2)$$
 is also an isomorphism.

 We have the following commutative exact diagram.
 \[\xymatrix{
  &  M\ar@{=}[d]  \ar[r] & (R_1\otimes_R M)\oplus (R_2\otimes_R M) \ar[d] \ar[r] & R'\otimes_{R_2} M \ar[d]^-{\delta}\ar[r] & 0\\
 0\ar[r] & M \ar[r] & X_1\oplus X_2 \ar[r] & R'\otimes_{R_2} X_2\ar[r] & 0
 }\]
 Here, the upper and lowers rows are given by (\ref{equ:M}) and the above pullback diagram, respectively; the middle vertical morphism, which is a diagonal matrix,  gives rise to the counit $\epsilon_{(X_1, X_2; c)}$. Since $\delta$ is an isomorphism, we infer the required statement by the snake lemma.
\end{proof}

\begin{prop}\label{prop:ext-Mil}
Assume that  $\pi_1\colon R_1\rightarrow R'$ is surjective and that ${\rm Ker}\pi_1$ is a universally superfluous ideal of $R_1$. Then the adjoint pair $({\rm Ind}, {\rm Pb})$ restricts to mutually inverse equivalences
$$\{M\in R\mbox{-{\rm Mod}}\; |\; M \mbox{ separated } \} \stackrel{\sim}\longrightarrow {\rm Tr}_{\rm gl}(\pi_1, \pi_2).$$
\end{prop}

\begin{proof}
Lemmas~\ref{lem:Au} and \ref{lem:Ba} actually describe the Auslander class and the Bass class of the adjoint pair $({\rm Ind}, {\rm Pb})$, respectively. Then the equivalences follow from the general standard fact on adjoint pairs \cite[Theorem~1.1]{FJ}.
\end{proof}

\begin{rem}
The above equivalences extend the ones in Theorem \ref{thm:Mil} under the additional condition that the ideal ${\rm Ker}\pi_1$  is universally superfluous. We do not know whether this condition is necessary.
\end{rem}

\section{A triangular matrix ring and a tilting module}
In this section, we associate an upper triangular matrix ring $\Gamma$  to the pullback diagram (\ref{equ:sq}) of rings. An explicit $1$-tilting $\Gamma$-module is given; see Proposition~\ref{prop:tilt}.

\subsection{A triangular matrix ring} We fix the pullback diagram (\ref{equ:sq}) of rings. We view $R'$ as an $R_1$-$R_2$-bimodule. Then its \emph{right $R_2$-dual}
$$R'^*={\rm Hom}_{R_2^{\rm op}}(R', R_2)$$
is naturally an $R_2$-$R_1$-bimodule. Here, the notation ${\rm Hom}_{R_2^{\rm op}}(-, -)$ means the Hom group taken in the category of right $R_2$-modules.  The evaluation map
\begin{align}\label{equ:ev}
{\rm ev}\colon R'^*\otimes_{R_1} R'\longrightarrow R_2, \quad f\otimes x'\mapsto f(x').
\end{align}
is an $R_2$-bimodule morphism.

 Consider the following upper triangular matrix ring
$$\Gamma=\begin{pmatrix} R_2 & R'^*\\
                            0 & R_1  \end{pmatrix}.$$
We will identify a  left $\Gamma$-module as a column vector $\begin{pmatrix} X_2 \\ X_1\end{pmatrix}$ equipped with the \emph{structure morphism} $\phi\colon R'^*\otimes_{R_1} X_1\rightarrow X_2$, where each $X_i$ is a left $R_i$-module and $\phi$   is an $R_2$-morphism.  The $\Gamma$-action is given as follows:
$$\begin{pmatrix} a_2 & f \\
                 0 & a_1 \end{pmatrix} \begin{pmatrix} x_2 \\ x_1\end{pmatrix}=\begin{pmatrix} a_2x_2+\phi(f\otimes x_1) \\ a_1x_1 \end{pmatrix}.$$
We suppress the structure morphism $\phi$ when it is clearly understood. For details, we refer to \cite[III.2]{ARS}.

For each triple $(X_1, X_2; c)$ in ${\rm Tr}(\pi_1, \pi_2)$, we consider the following composite
\begin{align}\label{equ:phi}
\phi\colon R'^*\otimes_{R_1}X_1 \longrightarrow R'^*\otimes_{R_1} (R'\otimes_{R_2} X_2) \xrightarrow{{\rm ev}\otimes_{R_2} X_2} R_2\otimes_{R_2} X_2 \stackrel{\sim}\longrightarrow X_2,
\end{align}
where the leftmost morphism sends $f\otimes x_1$ to $f\otimes c(1\otimes x_1)$, and the rightmost morphism is the canonical isomorphism. More precisely, if $c(1\otimes x_1)=\sum a'_i\otimes y_i\in R'\otimes_{R_2} X_2$, we have
$$\phi(f\otimes x_1)=\sum f(a'_i)y_i.$$
Then we obtain a left $\Gamma$-module $\begin{pmatrix} X_2 \\ X_1\end{pmatrix}$ with its  structure morphism $\phi$. This defines a functor
$$\Phi\colon {\rm Tr}(\pi_1, \pi_2) \longrightarrow \Gamma \mbox{-{\rm Mod}}, \quad (X_1, X_2; c)\mapsto \begin{pmatrix} X_2 \\ X_1\end{pmatrix}.$$
We observe that $\Phi$ is faithful and additive.

\begin{prop}\label{prop:phi}
Assume that $R'$ is finitely generated projective as the induced right $R_2$-module by $\pi_2$. Then the above functor $\Phi$ is an equivalence.
\end{prop}

\begin{proof}
By the assumption, we have a natural isomorphism of left $R'$-modules
$$R'\otimes_{R_2} X_2\longrightarrow {\rm Hom}_{R_2}(R'^*, X_2), \quad a'\otimes x_2\mapsto (f\mapsto f(a')x_2).$$
We infer that  there is a composition of natural isomorphisms
\begin{align*}
{\rm Hom}_{R'}(R'\otimes_{R_1}X_1, R'\otimes_{R_2}X_2)&\simeq {\rm Hom}_{R_1}(X_1, R'\otimes_{R_2}X_2)\\
& \simeq  {\rm Hom}_{R_1}(X_1, {\rm Hom}_{R_2}(R'^*, X_2))\\
&\simeq  {\rm Hom}_{R_2}(R'^*\otimes_{R_1} X_1, X_2),
\end{align*}
which sends the above $c$ to $\phi$. Then the verification of the equivalence is routine.
\end{proof}

Consider the following composite
$$R\mbox{-{\rm Mod}} \stackrel{{\rm Ind}}\longrightarrow {\rm Tr}(\pi_1, \pi_2) \stackrel{\Phi}\longrightarrow \Gamma\mbox{-{\rm Mod}}.$$
It sends $R$ to the $\Gamma$-module $\begin{pmatrix} R_2 \\ R_1\end{pmatrix}$, whose structure morphism is given by
\begin{align}\label{equ:psi}
\psi \colon R'^*\otimes_{R_1} R_1\longrightarrow R_2, \quad f\otimes 1 \rightarrow f(1).
\end{align}
Since there is a natural $R$-action on  $\begin{pmatrix} R_2 \\ R_1\end{pmatrix}$ from the right side, then $\begin{pmatrix} R_2 \\ R_1\end{pmatrix}$ is a naturally $\Gamma$-$R$-bimodule.

\begin{lem}
The composite $\Phi\circ  {\rm Ind}\colon R\mbox{-{\rm Mod}}\rightarrow \Gamma\mbox{-{\rm Mod}}$ is naturally isomorphic to the tensor functor $\begin{pmatrix} R_2 \\ R_1\end{pmatrix}\otimes_R-$.
\end{lem}

\begin{proof}
It suffices to observe that the composite sends a left $R$-module $M$ to the left $\Gamma$-module $\begin{pmatrix} R_2\otimes_R M \\ R_1\otimes_R M\end{pmatrix}$ with the structure morphism
$$\psi\otimes_R M\colon R'^*\otimes_{R_1} (R_1\otimes_R, M)\longrightarrow  R_2\otimes_R M,$$
where $\psi$ is given in (\ref{equ:psi}).
\end{proof}

\begin{rem}\label{rem:R}
Assume that $R'$ is finitely generated projective as the induced right $R_2$-module by $\pi_2$. In view of Remark~\ref{rem:Mil}(1) and Proposition~\ref{prop:phi}, we infer that  $\Phi\circ  {\rm Ind}$ restricts to a fully faithful functor on $R\mbox{-Proj}$. It follows that the  $\Gamma$-$R$-bimodule $\begin{pmatrix} R_2 \\ R_1\end{pmatrix}$ induces an isomorphism of rings $R^{\rm op}\simeq {\rm End}_\Gamma (\begin{pmatrix} R_2 \\ R_1\end{pmatrix}).$
\end{rem}

\subsection{An explicit $1$-tilting module}

Let us recall the notion of a tilting module.  For tilting modules and derived categories, we refer to \cite[Chapter 8]{KZ} and \cite[Chapters 4 and 5]{AHK}.

Let $\Lambda$ be a ring. For a left $\Lambda$-module $M$, we denote by ${\rm pd}_\Lambda(M)$ its projective dimension, and by ${\rm add}M$ the full subcategory of $\Lambda\mbox{-Mod}$ formed by direct summands of finite direct sums of $M$.

A finitely presented $\Lambda$-module $T$ is \emph{partial $1$-tilting}, provided that ${\rm pd}_\Lambda(T)\leq 1$ and ${\rm Ext}_\Lambda^1(T, T)=0$. A partial $1$-tilting $\Lambda$-module $T$ is \emph{$1$-tilting} \cite{HR},  if in addition there is an exact sequence of $\Lambda$-modules
$$0\longrightarrow \Lambda \longrightarrow  T^0 \longrightarrow T^1 \longrightarrow 0$$
with each $T^i\in {\rm add}T$.

\begin{rem}\label{rem:tilt}
Denote by $\langle T\rangle$ the smallest full subcategory of $\Lambda\mbox{-Mod}$ which contains $T$ and is closed under direct summands, extensions and kernels of epimorphisms. Then a partial $1$-tilting $\Lambda$-module $T$ is $1$-tilting if and only if $\Lambda\in \langle T\rangle$.

Indeed, the ``only if" part is clear. Conversely, if $\Lambda\in \langle T\rangle$,  then $T$, viewed as a stalk complex concentrated in degree zero, is a tilting complex over $\Lambda$.  Then the ``if" part follows from the well-known fact:  any module is  a $1$-tilting module if and only if it has projective dimension at most one and is a tilting complex.
\end{rem}

The following result is a very special case of \cite[Theorem 8.3.3]{KZ}. We denote by $\mathbf{D}(\Lambda\mbox{-Mod})$ the unbounded derived category of $\Lambda\mbox{-Mod}$. As usual, a complex of $\Lambda$-modules is denoted by $X^\bullet=(X^i, d_X^i)_{i \in \mathbb{Z}}$.

\begin{lem}\label{lem:ff}
Let $T$ be a partial $1$-tilting $\Lambda$-module. Consider $T$ as a $\Lambda$-$S$-bimodule with $S={\rm End}_\Lambda(T)^{\rm op}$. Then the derived tensor functor
$$T\otimes_S^\mathbb{L}- \colon \mathbf{D}(S\mbox{-{\rm Mod}})\longrightarrow \mathbf{D}(\Lambda\mbox{-{\rm Mod}})$$
is fully faithful. Moreover, if $T$ is $1$-tilting, then $T\otimes_S^\mathbb{L}-$ is an equivalence. \hfill $\square$
\end{lem}

Recall the left $\Gamma$-module $\begin{pmatrix} R_2 \\ R_1\end{pmatrix}$ from the previous subsection. We have another $\Gamma$-module $\begin{pmatrix} 0 \\ R_1\end{pmatrix}$ with the trivial structure morphism.

\begin{prop}\label{prop:tilt}
Assume that $\pi_1\colon R_1\rightarrow R'$ is surjective and that $R'$ is finitely generated projective as the induced right $R_2$-module by $\pi_2$. Then the $\Gamma$-module $\begin{pmatrix} R_2 \\ R_1\end{pmatrix}\oplus \begin{pmatrix} 0 \\ R_1\end{pmatrix}$ is $1$-tilting.
\end{prop}

\begin{proof}
 Set $T=\begin{pmatrix} R_2 \\ R_1\end{pmatrix}\oplus \begin{pmatrix} 0 \\ R_1\end{pmatrix}$. From the assumption, we infer that the $R_2$-module $R'^*={\rm Hom}_{R_2^{\rm op}}(R', R_2)$ is finitely generated projective. It follows that  the $\Gamma$-module $\begin{pmatrix} R'^* \\ 0\end{pmatrix}$ is also finitely generated projective.
 The following natural exact sequence of $\Gamma$-modules
\begin{align}\label{equ:exact1}
0 \longrightarrow \begin{pmatrix} R'^* \\0 \end{pmatrix} \longrightarrow \begin{pmatrix} R'^* \\ R_1\end{pmatrix} \longrightarrow \begin{pmatrix} 0 \\ R_1\end{pmatrix}\longrightarrow 0
\end{align}
implies that $\begin{pmatrix} 0 \\ R_1\end{pmatrix}$ is finitely presented with projective dimension at most one. Similarly, we have the following exact sequence
\begin{align}\label{equ:exact2}
0\longrightarrow \begin{pmatrix} R_2 \\ 0\end{pmatrix} \longrightarrow \begin{pmatrix} R_2 \\ R_1\end{pmatrix} \longrightarrow \begin{pmatrix} 0 \\ R_1\end{pmatrix} \longrightarrow 0.
\end{align}
It follows that $\begin{pmatrix} R_2 \\ R_1\end{pmatrix}$ is also finitely presented with projective dimension at most one.  This proves that $T$ is finitely presented satisfying ${\rm pd}_\Gamma(T)\leq 1$.

Since ${\rm Hom}_\Gamma( \begin{pmatrix} R'^* \\0 \end{pmatrix},  \begin{pmatrix} 0\\R_1 \end{pmatrix})=0$, we infer from (\ref{equ:exact1}) that ${\rm Ext}^1_\Gamma(\begin{pmatrix} 0\\R_1 \end{pmatrix}, \begin{pmatrix} 0\\R_1 \end{pmatrix})=0$. By the surjectivity of $\pi_1$ and the natural isomorphism $R'\simeq {\rm Hom}_{R_2}(R'^*, R_2)$, we infer that the following map
$$\iota\colon R_1\longrightarrow {\rm Hom}_{R_2}(R'^*, R_2), \quad a_1\mapsto (f\mapsto f\circ \pi_1(a_1))$$
is surjective. Take a $\Gamma$-morphism $\begin{pmatrix} g \\ 0\end{pmatrix}\colon \begin{pmatrix} R'^* \\ 0\end{pmatrix}\rightarrow \begin{pmatrix} R_2 \\ R_1\end{pmatrix}$, where $g\colon R'^*\rightarrow R_2$ is an $R_2$-morphism. Take $a_1\in R_1$ such that $\iota(a_1)=g$. Denote by $r_{a_1}\colon R_1\rightarrow R_1$ the map sending  $x$ to $xa_1$.  Then the following $\Gamma$-morphism
$$\begin{pmatrix} g \\ r_{a_1}\end{pmatrix}\colon \begin{pmatrix} R'^* \\ R_1\end{pmatrix} \longrightarrow \begin{pmatrix} R_2 \\ R_1\end{pmatrix}$$
is well defined, and proves that the given $\Gamma$-morphism $\begin{pmatrix} g \\ 0\end{pmatrix}$ factors through the inclusion in (\ref{equ:exact1}). It follows that ${\rm Ext}^1_\Gamma(\begin{pmatrix} 0 \\ R_1\end{pmatrix}, \begin{pmatrix} R_2 \\ R_1\end{pmatrix})=0$. Therefore, we have ${\rm Ext}_\Gamma^1(\begin{pmatrix} 0 \\ R_1\end{pmatrix}, T)=0$. Apply ${\rm Ext}_\Gamma^1(-, T)$ to (\ref{equ:exact2}), we infer that ${\rm Ext}_\Gamma^1(\begin{pmatrix} R_2 \\ R_1\end{pmatrix}, T)=0$, proving ${\rm Ext}^1_\Gamma(T,T)=0$.

We observe from (\ref{equ:exact1}) and (\ref{equ:exact2}) that $\Gamma\simeq \begin{pmatrix} R_2 \\ 0\end{pmatrix}\oplus \begin{pmatrix} R'^* \\ R_1\end{pmatrix} $ lies in $\langle T\rangle$, the smallest full subcategory of $\Gamma\mbox{-Mod}$ containing $T$ and closed under direct summands, extensions and kernels of epimorphisms. Then Remark \ref{rem:tilt} implies that $T$ is $1$-tilting.
\end{proof}

Keep the assumptions of Proposition~\ref{prop:tilt}. Then the opposite  ring  of the endomorphism ring of $\begin{pmatrix} R_2 \\ R_1\end{pmatrix}\oplus \begin{pmatrix} 0 \\ R_1\end{pmatrix}$ is isomorphic to the following matrix ring.
$$\Gamma'=\begin{pmatrix} R & R_1\\
                          I_1 & R_1 \end{pmatrix}$$
Here, $I_1={\rm Ker}\pi_1$. The corresponding morphism $I_1\otimes_{R_1} R_1\rightarrow R_1$ is the canonical embedding. For the other morphism $R_1\otimes_R I_1\rightarrow R$, we use the isomorphism between $I_1$ and  ${\rm Ker}i_2$, which is restricted from $i_1$. We omit the details of the computation. We mention that \cite[Section~3]{BD} studies a matrix algebra of the same form as $\Gamma'$.

The following consequence is immediate from Proposition~\ref{prop:tilt} and Lemma~\ref{lem:ff}.

\begin{cor}
Keep the assumptions in Proposition \ref{prop:tilt}. Then the above $1$-tilting $\Gamma$-module induces a derived equivalence
$$\mathbf{D}(\Gamma \mbox{-{\rm Mod}})\simeq \mathbf{D}(\Gamma' \mbox{-{\rm Mod}}).$$
\end{cor}

By \cite[Corollary~4.11]{Lad},  the ring $\Gamma$ is derived equivalent to another upper triangular matrix ring
$$\Gamma''=\begin{pmatrix} R_1 & R'\\
                         0 & R_2 \end{pmatrix}$$
via a two-term tilting complex over $\Gamma$,  where $R'$ is again viewed as an $R_1$-$R_2$-bimodule. Therefore, the three matrix rings $\Gamma$, $\Gamma'$ and $\Gamma''$ are all derived equivalent.

In the next section, we will need the following result, which is also an immediate consequence of Proposition~\ref{prop:tilt} and Lemma~\ref{lem:ff}.

\begin{cor}\label{cor:ff}
Keep the assumptions in Proposition \ref{prop:tilt}. Then the following derived tensor functor
$$\begin{pmatrix}R_2\\ R_1\end{pmatrix}\otimes^\mathbb{L}_R- \colon \mathbf{D}(R \mbox{-{\rm Mod}})\longrightarrow \mathbf{D}(\Gamma \mbox{-{\rm Mod}})$$
is fully faithful.
\end{cor}

\begin{proof}
We use the fact that  the natural ring homomorphism $R\rightarrow {\rm End}_\Gamma (\begin{pmatrix}R_2\\ R_1\end{pmatrix})^{\rm op}$,  induced by the $\Gamma$-$R$-bimodule $\begin{pmatrix}R_2\\ R_1\end{pmatrix}$,  is an isomorphism; see Remark~\ref{rem:R}.
\end{proof}

\section{The main result and its proof}

In this section, we prove the main result, which establishes an epivalence from the derived category of the pullback ring to a certain subcategory of gluing derived-triples; see Theorem~\ref{thm:main}. The proof relies on Corollary~\ref{cor:ff} and a general result on recollements in \cite{CL}.

\subsection{The comma category and derived-triples}

Let $F\colon \mathcal{C}\rightarrow \mathcal{D}$ be a functor. The functor $F$ is called an \emph{epivalence},  provided that it is full and dense, and detects isomorphisms. The last condition means that any morphism $f$ in $\mathcal{C}$ is an isomorphism if and only if $F(f)$ is an isomorphism in $\mathcal{D}$.

Assume that both $\mathcal{C}$ and $\mathcal{D}$ are additive categories and that $F$ is an additive functor. The \emph{kernel ideal} of $F$ means the two-sided ideal of $\mathcal{C}$ formed by those morphisms annihilated by $F$.

The following fact is standard.

\begin{lem}\label{lem:epiv}
Let  $F\colon \mathcal{C}\rightarrow \mathcal{D}$ be an additive functor. Assume that $F$ is full such that its kernel ideal is nilpotent. Then $F$ induces an epivalence $\mathcal{C}\rightarrow {\rm Im}(F)$, where ${\rm Im}(F)$ denotes the essential image of $F$.
\end{lem}

\begin{proof}
It suffices to prove the following statement:  any morphism $f\colon X\rightarrow Y$ in $\mathcal{C}$ is an isomorphism provided that $F(f)$ is an isomorphism in $\mathcal{D}$. By the fullness of $F$, there exists a morphism $g\colon Y\rightarrow X$ such that $F(g\circ f)={\rm Id}_{F(X)}$. It follows that ${\rm Id}_X-g\circ f$ lies in the kernel ideal of $F$.  By the nilpotent property of the kernel ideal, we infer that $g\circ f$ is  an isomorphism, that is, $f$ has a left inverse. Similarly, $f$ has a right inverse. Therefore, the morphism $f$ is an isomorphism, as required.
\end{proof}

For a functor $F\colon \mathcal{C}\rightarrow \mathcal{D}$, the \emph{comma category} $(F\downarrow \mathcal{D})$ is defined as follows:  its objects are triples $(C, D; f)$, where $C$ and $D$ are objects in $\mathcal{C}$ and $\mathcal{D}$, respectively, and $f\colon F(C)\rightarrow D$ is a morphism in $\mathcal{D}$; a morphism $(u, v)\colon (C, D; f)\rightarrow (C', D'; f')$ consists of a morphism $u\colon C\rightarrow C'$ in $\mathcal{C}$ and a morphism $v\colon D\rightarrow D'$ in $\mathcal{D}$ satisfying $f'\circ F(u)=v\circ f$.

We recall the pullback diagram (\ref{equ:sq}) and the upper triangular matrix ring $\Gamma$ from Section~3.  We will consider the derived tensor functor
$$R'^*\otimes^\mathbb{L}_{R_1}-\colon \mathbf{D}(R_1\mbox{-Mod})\longrightarrow \mathbf{D}(R_2\mbox{-Mod})$$
and its comma category $( R'^*\otimes^\mathbb{L}_{R_1}- \downarrow \mathbf{D}(R_2\mbox{-Mod}))$.

Recall that a left $\Gamma$-module is viewed as a column vector $\begin{pmatrix} X_2\\ X_1\end{pmatrix}$ with a structure $R_2$-morphism $\phi\colon R'^*\otimes_{R_1} X_1\rightarrow X_2$. This shows that the comma category of the functor $R'^*\otimes_{R_1}-\colon R_1\mbox{-Mod}\rightarrow R_2\mbox{-Mod}$ is equivalent to $\Gamma\mbox{-Mod}$.

Similarly, a complex of $\Gamma$-modules is viewed as a column vector  $\begin{pmatrix} X^\bullet_2\\ X^\bullet_1\end{pmatrix}$ with a chain map $\phi^\bullet\colon R'^*\otimes_{R_1} X^\bullet_1\rightarrow X^\bullet_2$ between complexes of $R_2$-modules, where each $X_i^\bullet$ is a complex of $R_i$-modules. Then we have a well-defined functor
$$\Psi\colon \mathbf{D}(\Gamma\mbox{-Mod})\longrightarrow ( R'^*\otimes^\mathbb{L}_{R_1}- \downarrow \mathbf{D}(R_2\mbox{-Mod})), \quad \begin{pmatrix} X^\bullet_2\\ X^\bullet_1\end{pmatrix}\mapsto (X_1^\bullet, X_2^\bullet; \tilde{\phi}^\bullet),$$
where $\tilde{\phi}^\bullet$ is the composition of $\phi^\bullet$ and the canonical morphism $R'^*\otimes^\mathbb{L}_{R_1}X^\bullet_1\rightarrow R'^*\otimes_{R_1} X^\bullet_1$.

The following result extends \cite[Appendix, Proposition~1]{Gei}.

\begin{prop}\label{prop:psi}
Keep the notation as above. Then the functor $\Psi$ is full and dense, whose kernel ideal is square zero. In particular, $\Psi$ is an epivalence.
\end{prop}

\begin{proof}
For the triangular matrix ring $\Gamma$, there is a canonical recollement \cite{BBD} as follows:
\[
\xymatrix{
\mathbf{D}(R_2\mbox{-{\rm Mod}}) \ar[rr]|{i} &&  \mathbf{D}(\Gamma\mbox{-{\rm Mod}})\ar[rr]|{j}  \ar@/_1pc/[ll]|{i_\lambda} \ar@/^1pc/[ll]|{i_\rho} &&  \mathbf{D}(R_1\mbox{-{\rm Mod}}).   \ar@/_1pc/[ll]|{j_\lambda} \ar@/^1pc/[ll]|{j_\rho}
}\]
The functors are given as follows:  $i(X_2^\bullet)=\begin{pmatrix} X^\bullet_2\\ 0\end{pmatrix}$, $j(\begin{pmatrix} X^\bullet_2\\ X^\bullet_1\end{pmatrix})=X_1^\bullet$, $i_\rho(\begin{pmatrix} X^\bullet_2\\ X^\bullet_1\end{pmatrix})=X_2^\bullet$, $j_\rho (X_1^\bullet)=\begin{pmatrix} 0\\ X^\bullet_1\end{pmatrix}$, and $j_\lambda(X_1^\bullet)=\begin{pmatrix} R'^*\otimes^\mathbb{L}_{R_1}X^\bullet_1\\ X^\bullet_1\end{pmatrix}$. Since there is a surjective ring homomorphism $\Gamma\rightarrow R_2$, we view $R_2$ as an  $R_2$-$\Gamma$-bimodule. Then $i_\lambda=R_2\otimes^\mathbb{L}_\Gamma-$.

We observe that the composition $i_\rho j_\lambda$ is isomorphic to $R'^*\otimes^\mathbb{L}_{R_1}-$. Now the result follows from a general result on recollements \cite[Theorem A]{CL} and Lemma~\ref{lem:epiv}.
\end{proof}

We will describe a derived analogue of Proposition~\ref{prop:phi}. A \emph{derived-triple} $(X_1^\bullet, X_2^\bullet; c^\bullet)$ consists of complexes $X_i^\bullet$ of $R_i$-modules and a morphism
$$c^\bullet \colon R'\otimes^\mathbb{L}_{R_1} X_1^\bullet \longrightarrow R'\otimes^\mathbb{L}_{R_2} X_2^\bullet$$ in $\mathbf{D}(R'\mbox{-Mod})$. The derived-triple is said to be \emph{gluing},  provided that $c^\bullet$ is an isomorphism in $\mathbf{D}(R'\mbox{-Mod})$. A morphism
$$(f_1^\bullet, f_2^\bullet)\colon (X_1^\bullet, X_2^\bullet; c^\bullet) \longrightarrow (X_1'^\bullet, X_2'^\bullet; {c'}^\bullet)$$
between derived-triples consists of morphisms $f_i^\bullet$ in $\mathbf{D}(R_i\mbox{-Mod})$ satisfying
$${c'}^\bullet \circ (R'\otimes^{\mathbb{L}}_{R_1} f_1)=(R'\otimes^{\mathbb{L}}_{R_2} f_2)\circ c^\bullet.$$
This defines the category  $\mathbf{D}{\rm Tr}(\pi_1, \pi_2)$ of  derived-triples.  Denote by $\mathbf{D}{\rm Tr}_{\rm gl}(\pi_1, \pi_2)$ the full subcategory formed by gluing derived-triples.

Assume that we are given a derived-triple $(X_1^\bullet, X_2^\bullet; c^\bullet)$. Similar to (\ref{equ:phi}), the following composite defines
a morphism $\phi^\bullet\colon R'^*\otimes^\mathbb{L}_{R_1} X_1^\bullet\rightarrow X_2^\bullet$ in $\mathbf{D}(R_2\mbox{-Mod})$.
\[
\xymatrix{
R'^* \otimes_{R_1}^\mathbb{L} X^\bullet_1 \ar[rr]^-{R'^* \otimes_{R_1}^\mathbb{L} {\rm can}}  &&
R'^*\otimes_{R_1}^\mathbb{L} (R'\otimes_{R_1}^\mathbb{L} X^\bullet_1) \ar[rr]^{R'^*\otimes_{R_1}^\mathbb{L} c^\bullet} &&
R'^*\otimes_{R_1}^\mathbb{L} (R'\otimes^{\mathbb{L}}_{R_2} X^\bullet_2) \ar[d]^-{\rm can} \\
 X^\bullet_2 && R_2\otimes_{R_2} X^\bullet_2 \ar[ll]_-{\rm can}&&  \ar[ll]_-{{\rm ev}\otimes_{R_2} X^\bullet_2} R'^*\otimes_{R_1} R'\otimes_{R_2} X^\bullet_2
    }\]
Here, the above three ``can" mean the canonical morphisms, and ``ev" is the evaluation map (\ref{equ:ev}). This yields a well-defined functor
$$\mathbf{D}\Phi\colon \mathbf{D}{\rm Tr}(\pi_1, \pi_2) \longrightarrow ( R'^*\otimes^\mathbb{L}_{R_1}- \downarrow \mathbf{D}(R_2\mbox{-Mod})),
$$
sending the derived-triple $(X_1^\bullet, X_2^\bullet; c^\bullet)$ to the corresponding object $(X_1^\bullet, X_2^\bullet; \phi^\bullet)$ in the comma category.

We omit the reasoning of the following result, since it is the same as the one in the proof of Proposition~\ref{prop:phi}.

\begin{prop}\label{prop:phi-der}
Assume that $R'$ is finitely generated projective as the induced right $R_2$-module by $\pi_2$. Then the above functor $\mathbf{D}\Phi$ is an equivalence. \hfill $\square$
\end{prop}

\subsection{The derived induction functor}
For each complex $M^\bullet$ of $R$-modules, we have a gluing derived-triple
$${\rm Ind}^\mathbb{L}(M^\bullet)=(R_1\otimes_{R}^\mathbb{L} M^\bullet,R_2\otimes_{R}^\mathbb{L} M^\bullet; {\rm can}_{M^\bullet}),$$
where ${\rm can}_{M^\bullet}$ is the canonical isomorphism. Therefore, we have the \emph{derived induction} functor
$${\rm Ind}^\mathbb{L}\colon \mathbf{D}(R\mbox{-Mod})\longrightarrow \mathbf{D}{\rm Tr}(\pi_1, \pi_2).$$

\begin{thm}\label{thm:main}
Assume that $\pi_1\colon R_1\rightarrow R'$ is surjective and that $R'$ is finitely generated projective as the induced right $R_2$-module by $\pi_2$. Then the derived induction functor ${\rm Ind}^\mathbb{L}$ is full, whose kernel ideal is square zero. In particular, it restricts to an epivalence
$$ \mathbf{D}(R\mbox{-{\rm Mod}})\longrightarrow {\rm Im}({\rm Ind}^\mathbb{L}),$$
where the essential image ${\rm Im}({\rm Ind}^\mathbb{L})$ of ${\rm Ind}^\mathbb{L}$ lies in $\mathbf{D}{\rm Tr}_{\rm gl}(\pi_1, \pi_2)$.
\end{thm}

\begin{proof}
We observe that the following diagram commutes up to a natural isomorphism.
\[\xymatrix{
\mathbf{D}(R\mbox{-{\rm Mod}}) \ar[d]_-{\begin{pmatrix}R_2\\ R_1\end{pmatrix}\otimes_R^\mathbb{L}-}  \ar[rr]^-{{\rm Ind}^\mathbb{L}} && \mathbf{D}{\rm Tr}(\pi_1, \pi_2) \ar[d]^-{\mathbf{D}\Phi}\\
\mathbf{D}(\Gamma\mbox{-{\rm Mod}}) \ar[rr]^-{\Psi} &&  ( R'^*\otimes^\mathbb{L}_{R_1}- \downarrow \mathbf{D}(R_2\mbox{-Mod}))
}\]
By Corollary~\ref{cor:ff}, the functor $\begin{pmatrix}R_2\\ R_1\end{pmatrix}\otimes_R^\mathbb{L}-$ is fully faithful. By Proposition~\ref{prop:psi}, the functor $\Psi$ is full and dense with a square-zero kernel ideal. The functor $\mathbf{D}\Phi$ is an equivalence by Proposition~\ref{prop:phi-der}. Then the required result follows immediately from the commutative square above and Lemma~\ref{lem:epiv}.
\end{proof}

In general, the essential image ${\rm Im}({\rm Ind}^\mathbb{L})$ is not easy to characterize. We now restricts to the derived category $\mathbf{D}^{-}(R\mbox{-Mod})$ of right-bounded complexes of $R$-modules, which is identified with the full subcategories of $\mathbf{D}(R\mbox{-Mod})$ formed by complexes with right-bounded cohomology. Similarly, we have the full subcategory $\mathbf{D}^{-}{\rm Tr}(\pi_1, \pi_2)$ and $\mathbf{D}^{-}{\rm Tr}_{\rm gl}(\pi_1, \pi_2)$ of $\mathbf{D}{\rm Tr}(\pi_1, \pi_2)$ given by right-bounded complexes.

Let $\Lambda$ be a left perfect ring. Then any $\Lambda$-module has a projective cover. Consequently, any complex $P^\bullet$ of projective modules is homotopically equivalent to a \emph{minimal} complex $Q^\bullet$ of projective modules. Here, the minimality means that the image of each differential $d_Q^n\colon Q^n\rightarrow Q^{n+1}$ lies in ${\rm rad}(Q^{n+1})={\rm rad}(\Lambda)Q^{n+1}$. For details, we refer to \cite[Appendix B]{Kra}.

\begin{prop}\label{prop:rb}
Assume that $\pi_1$ is surjective and that $R'$ is finitely generated projective as the induced right $R_2$-module by $\pi_2$.  Assume further that both $R_1$ and $R_2$ are left perfect such that $\pi_2({\rm rad}(R_2))\subseteq {\rm rad}(R')$. Then the derived induction functor
$${\rm Ind}^\mathbb{L}\colon \mathbf{D}^{-}(R\mbox{-{\rm Mod}})\longrightarrow \mathbf{D}^{-}{\rm Tr}_{\rm gl}(\pi_1, \pi_2)$$
is full and dense, whose kernel ideal is square zero. In particular, it is an epivalence.
\end{prop}

\begin{proof}
In view of Theorem~\ref{thm:main}, it suffices to prove the density. Take any derived-triple $(X_1^\bullet, X_2^\bullet; c^\bullet)$ in $\mathbf{D}^{-}{\rm Tr}_{\rm gl}(\pi_1, \pi_2)$. Taking projective resolutions and replacing the resolutions by some minimal complexes, we may assume that $X_i^\bullet=P_i^\bullet$ is a minimal right-bounded complex of projective $R_i$-modules.

 We recall the assumption $\pi_2({\rm rad}(R_2))\subseteq {\rm rad}(R')$ and observe $\pi_1({\rm rad}(R_1))\subseteq {\rm rad}(R')$ by the surjectivity of $\pi_1$. Consequently, both the complexes $R'\otimes^\mathbb{L}_{R_i} P_i^\bullet\simeq R'\otimes_{R_i}P^\bullet_i$ of projective $R'$-modules are minimal. Therefore, the isomorphism $c^\bullet$ is represented by a chain isomorphism over $R'$. Without loss of generality, we assume that
 $$c^\bullet\colon R'\otimes_{R_1} P^\bullet_1\longrightarrow R'\otimes_{R_2} P^\bullet_2$$
 is a chain isomorphism. Hence, on each degree $i$,
 $$c^i\colon R'\otimes_{R_1} P_1^i \longrightarrow R'\otimes_{R_2} P_2^i$$
  is an isomorphism of $R'$-modules. We apply Theorem~\ref{thm:Mil} to obtain a projective $R$-module $P^i$ such that ${\rm Ind}(P^i)\simeq (P_1^i, P_2^i; c^i)$. These projective $R$-modules $P^i$ form a right-bounded complex $P^\bullet$ of projective $R$-modules. We observe that ${\rm Ind}^\mathbb{L}(P^\bullet)$ is isomorphic to $(P_1^\bullet, P_2^\bullet; c^\bullet)=(X_1^\bullet, X_2^\bullet; c^\bullet)$, completing the proof.
\end{proof}

\begin{rem}\label{rem:image}
(1) Keep the same assumptions as in Proposition~\ref{prop:rb}. Suppose that we are given a gluing derived-triple $(X_1^\bullet, X_2^\bullet; c^\bullet)$ of unbounded complexes. We can still obtain a complex $P^\bullet$ of projective $R$-modules, which is unbounded and thus might be not homotopically projective \cite[Subsection 8.1.1]{KZ}. So, it is not clear whether $R_i\otimes_R P^\bullet$ is isomorphic to $R_i\otimes^\mathbb{L}_R P^\bullet$. However, if we suppose in addition that as right $R$-modules, both $R_1$ and $R_2$ have finite flat dimension, then these two complexes are indeed isomorphic. Consequently, we still have an isomorphism ${\rm Ind}^\mathbb{L}(P^\bullet)\simeq (X_1^\bullet, X_2^\bullet; c^\bullet)$.

In conclusion, if in addition both $R_1$ and $R_2$ have finite flat dimension as right $R$-modules, the derived induction functor
$${\rm Ind}^\mathbb{L}\colon \mathbf{D}(R\mbox{-{\rm Mod}})\longrightarrow \mathbf{D}{\rm Tr}_{\rm gl}(\pi_1, \pi_2)$$
is dense, and thus an epivalence.

(2) Let us consider the derived categories of right-bounded complexes of  finitely presented modules.  We replace the left perfectness condition of $R_i$ in Proposition~\ref{prop:rb} by a weaker  condition that each $R_i$ is semiperfect. Then by the same proof of Proposition~\ref{prop:rb}  and using Theorem~\ref{thm:Mil} on $R\mbox{-proj}$, we infer that the derived induction functor
$${\rm Ind}^\mathbb{L}\colon \mathbf{D}^{-}(R\mbox{-{\rm mod}})\longrightarrow \{(X_1^\bullet, X_2^\bullet; c^\bullet)\in \mathbf{D}^{-}{\rm Tr}_{\rm gl}(\pi_1, \pi_2)\; |\; X_i^\bullet\in \mathbf{D}^{-}(R_i\mbox{-{\rm mod}})\}$$
is full and dense, whose kernel ideal is square zero. Here, for each ring $\Lambda$ we identify $\mathbf{D}^{-}(\Lambda\mbox{-{\rm mod}})$ with  the homotopy category $\mathbf{K}^{-}(\Lambda\mbox{-proj})$ of right bounded complexes of finitely generated projective modules. This epivalence generalizes and strengthens  \cite[Theorem~2.4]{BD04}; compare \cite[Theorem~4.2]{BD}.
\end{rem}

\vskip 5pt

\noindent{\bf Acknowledgement.}\quad X.W. Chen is inspired by a talk of Igor Burban,  presented on the Oberwolfach workshop in January 2020. We are very grateful to Hongxing Chen, Chao Zhang and Yu Zhou for helpful discussion. This work is supported by National Natural Science Foundation of China (No.s 11671174, 11671245 and 11971449).

\bibliography{}

\vskip 10pt

 {\footnotesize \noindent Xiao-Wu Chen, Jue Le\\
 Key Laboratory of Wu Wen-Tsun Mathematics, Chinese Academy of Sciences,\\
 School of Mathematical Sciences, University of Science and Technology of China, Hefei 230026, Anhui, PR China}

\end{document}